\documentclass[11pt]{article}
\usepackage{amscd}
\usepackage{amsfonts}
\usepackage{amsmath}
\usepackage{amssymb}
\usepackage{amsthm}
\usepackage{bbm}
\usepackage{CJK}
\usepackage{fancyhdr}
\usepackage{graphicx}
\usepackage{hyperref}
\usepackage{indentfirst}
\usepackage{latexsym}
\usepackage{mathrsfs}
\usepackage{xypic}

\usepackage[top=1in,bottom=1in,left=1.25in,right=1.25in]{geometry}
\def\<{\langle}
\def\>{\rangle}
\def\a{\alpha}
\def\b{\beta}

\def\c{\cdot}

\def\g{\gamma}

\def\o{\otimes}

\date{}
\begin{document}
\renewcommand{\baselinestretch}{1.2}
\renewcommand{\arraystretch}{1.0}
\title{\bf Cohomology and derivations of  BiHom-Lie conformal algebras}
\author{{\bf Shuangjian Guo$^{1}$, Xiaohui Zhang$^{2}$,  Shengxiang Wang$^{3}$\footnote
        { Corresponding author(Shengxiang Wang):~~wangsx-math@163.com} }\\
{\small 1. School of Mathematics and Statistics, Guizhou University of Finance and Economics} \\
{\small  Guiyang  550025, P. R. of China} \\
{\small 2.  School of Mathematical Sciences, Qufu Normal University}\\
{\small Qufu  273165, P. R. of China}\\
{\small 3.~ School of Mathematics and Finance, Chuzhou University}\\
 {\small   Chuzhou 239000,  P. R. of China}}
 \maketitle
\begin{center}
\begin{minipage}{13.cm}

{\bf \begin{center} ABSTRACT \end{center}}
 In this paper,  we introduce the notion  of a BiHom-Lie conformal algebra and develop its cohomology theory.
Also,  we   discuss some applications to the study of deformations of regular BiHom-Lie conformal algebras.
 Finally, we introduce derivations of multiplicative
BiHom-Lie conformal algebras and study their properties.


 \end{minipage}
 \end{center}
 \normalsize\vskip1cm

\section*{INTRODUCTION}
\def\theequation{0. \arabic{equation}}
\setcounter{equation} {0}

Lie conformal algebras were introduced by Kac in \cite{Kac98},
he gave an axiomatic description of the singular part of the operator product expansion of chiral fields in conformal field theory.
It is an useful tool to study vertex algebras and has many applications in the theory of Lie algebras. Moreover.
Lie conformal algebras  have close connections to Hamiltonian formalism in the theory of nonlinear evolution equations. 
Lie conformal algebras are widely studied in the following aspects: The structure theory \cite{D'Andrea1998}, 
representation theory \cite{Cheng1997}, \cite{Cheng1998} and cohomology
theory \cite{Bakalov1999} of finite Lie conformal algebras have been developed. Moreover, Liberati
in \cite{L08} introduced a conformal analog of Lie bialgebras including the conformal classical Yang-Baxter equation, the conformal Manin triples and conformal Drinfeld¡¯s double.

 The notion of Hom-Lie algebras was firstly introduced by Hartwig,
Larsson and Silvestrov  to describe the structure of certain $q$-deformations
of the Witt and the Virasoro algebras, see \cite{Aizawa,  Hartwig, Hu},
which  was motivated by
applications to physics and to deformations of Lie algebras, especially Lie algebras
of vector fields.

  Recently, the Hom-Lie conformal algebra was introduced
and studied in \cite{Yuan14}, where it was proved that a Hom-Lie conformal algebra is equivalent to a Hom-Gel'fand-Dorfman bialgebra.  
Later, Zhao, Yuan and Chen  developed the cohomology theory of Hom-Lie conformal algebras and discuss
some applications to the study of deformations of regular Hom-Lie conformal
algebras. Also, they  introduced  $\a^k$-derivations of multiplicative Hom-Lie conformal
algebras and study their properties in \cite{Zhao2016}.

A BiHom-algebra is an algebra in such a way that the identities defining the structure
are twisted by two homomorphisms $\a,\b$. This class of algebras was introduced from a
categorical approach in \cite{Graziani} as an extension of the class of Hom-algebras. When the two
linear maps are same automorphisms, Bihom-algebras will be return to Hom-algebras.
These algebraic structures include BiHom-associative algebras, BiHom-Lie algebras and
BiHom-bialgebras.

Motivated by these results, the present paper is organized as follows. 
In Section 2,  we introduce the notion of a BiHom-Lie conformal algebra and  develop the cohomology theory of BiHom-Lie conformal algebras. 
In Section 3, we  discuss some applications to the study of deformations of regular BiHom-Lie conformal algebras. 
In Section 4, we study derivations of multiplicative BiHom-Lie conformal algebras and prove that the direct sum of the space of derivations is a BiHom-Lie
conformal algebra. In particular, any derivation gives rise to a derivation extension of a multiplicative BiHom-Lie conformal algebra. 
In Section 5, we introduce  generalized derivations of multiplicative BiHom-Lie conformal algebras  and study their properties.
\section{Preliminaries}
\def\theequation{\arabic{section}.\arabic{equation}}
\setcounter{equation} {0}

Throughout the paper, all algebraic systems are supposed to be over a field $\mathbb{C}$. 
We denote by $\mathbb{Z}$ the set of all integers and $\mathbb{Z}_{+}$ the set of all nonnegative integers.

In this section we recall some basic definitions and results related to our paper from \cite{Graziani} and \cite{Yuan14} .
\medskip

\noindent{\bf 1.1.  Hom-Lie conformal algebra} A Hom-Lie conformal algebra is a  $\mathbb{C}[\partial]$-module $R$ equipped with a linear endomorphism
$\a$ such that $\a\partial=\partial\a$, and a $\lambda$-bracket $[\c_{\lambda}\c]$ which defines a $\mathbb{C}$-bilinear map from $R\o R$ to  $R[\lambda]=\mathbb{C}[\lambda]\o R$ such that the following axioms hold:
\begin{eqnarray*}
&&[\partial a_{\lambda}b] =-\lambda[a_\lambda b],[ a_{\lambda}\partial b] =(\partial+\lambda)[a_\lambda b],\\
&&[a_{\lambda}b]=-[b_{-\lambda-\partial}a],\\
&&[\a(a)_\lambda[b_\mu c]]=[[a_{\lambda}b]_{\lambda+\mu}\a(c)]+[\a(b)_{\mu}[a_{\lambda}c]],
\end{eqnarray*}
 for $a,b,c\in R$.
\medskip

\noindent{\bf 1.2. BiHom-Lie algebra }  A BiHom-Lie algebra
is a 4-tuple $(L,[\cdot,\cdot],\alpha,\beta)$,
where $L$ is a ${k}$-linear space, $\alpha, \beta: L\rightarrow L$ and
 $[\cdot,\cdot]: L\o L\rightarrow L$ are linear maps,
 satisfying the following conditions:
\begin{eqnarray*}
&&\alpha\circ\beta=\beta\circ\alpha,\\
&&\alpha[a,a']=[\alpha(a),\alpha(a')],\beta[a,a']=[\beta(a),\beta(a')],\\
&&[\beta(a),\alpha(a')]=-[\beta(a'),\alpha(a)],\\
&&[\beta^{2}(a),[\beta(a'),\alpha(a'')]]+[\beta^{2}(a'),[\beta(a''),\alpha(a)]]+[\beta^{2}(a''),[\beta(a),\alpha(a')]]=0,
\end{eqnarray*}
 for all  $a,a',a''\in A$.
 \medskip
 
A BiHom-Lie  algebra $(L,[\cdot,\cdot],\alpha,\beta)$ is called multiplicative if $\a,\b$ are algebra endomorphisms,
i.e., $\a([a, b]) = [\a(a),  \a(b)], \b([a, b]) = [\b(a), \b(b)]$, for any $a,b\in L$. In particular, if $\a,\b$ are algebra
isomorphisms, then $(L,[\cdot,\cdot],\alpha,\beta)$ is called regular.
\medskip

Obviously, a Hom-Lie algebra $(L,[\cdot,\cdot],\alpha)$ is a particular case of a BiHom-Lie
algebra, namely $(L,[\cdot,\cdot],\alpha,\alpha)$.
 Conversely, a BiHom-Lie algebra $(L,[\cdot,\cdot],\alpha,\alpha)$ with bijective $\alpha$
is the Hom-Lie algebra $(L,[\cdot,\cdot],\alpha)$ .

\section{ Cohomology  of BiHom-Lie conformal algebras }
\def\theequation{\arabic{section}. \arabic{equation}}
\setcounter{equation} {0}

In this section, we introduce the notation  of BiHom-Lie conformal algebras and develop the cohomology theory of BiHom-Lie conformal algebras.
\medskip

\noindent{\bf Definition 2.1.} A BiHom-Lie conformal algebra is a  $\mathbb{C}[\partial]$-module $R$ equipped with two commuting  linear maps
$\a,\b$ such that $\a\partial=\partial\a, \b\partial=\partial\b, \a([a_\lambda b]) = [\a(a)_\lambda \a(b)], \b([a_\lambda b]) = [\b(a)_\lambda \b(b)]$, and a $\lambda$-bracket $[\c_{\lambda}\c]$ which defines a $\mathbb{C}$-bilinear map from $R\o R$ to  $R[\lambda]=\mathbb{C}[\lambda]\o R$ such that the following axioms hold for $a,b,c\in R$:
\begin{eqnarray}
&&[\partial a_{\lambda}b] =-\lambda[a_\lambda b], [ a_{\lambda}\partial b] =(\partial+\lambda)[a_\lambda b],\\
&&[\b(a)_{\lambda}\a(b)]=-[\b(b)_{-\lambda-\partial}\a(a)],\\
&&[\a\b(a)_\lambda[b_\mu c]]=[[\b(a)_{\lambda}b]_{\lambda+\mu} \b(c)]+[\b(b)_{\mu}[\a(a)_{\lambda}c]].
\end{eqnarray}

A BiHom-Lie conformal algebra $(R,\a,\b)$ is called multiplicative if $\a,\b$ are algebra endomorphisms,
i.e., $\a([a_\lambda b]) = [\a(a)_\lambda \a(b)], \b([a_\lambda b]) = [\b(a)_\lambda \b(b)]$, for any $a,b\in R$. 
In particular, if $\a,\b$ are algebra
isomorphisms, then $(R,\a,\b)$ is called regular.
\medskip

A BiHom-Lie conformal algebra $R$ is called finite if $R$ is a finitely generated $\mathbb{C}[\partial]$-module. The
rank of $R$ is its rank as a $\mathbb{C}[\partial]$-module.
\medskip

 \noindent{\bf Example 2.2.} Let $R$ be a Lie conformal algebra and $\a,\b:A \rightarrow A$ two commuting linear maps such that 
 $$\a\partial=\partial\a, \b\partial=\partial\b, \a([a_\lambda b]) = [\a(a)_\lambda \a(b)], \b([a_\lambda b]) = [\b(a)_\lambda \b(b)].$$ 
Define   a $\mathbb{C}$-bilinear map from $R\o R$ to  $R[\lambda]=\mathbb{C}[\lambda]\o R$ by $[a_{\lambda}b]'=[\a(a)_{\lambda}\b(b)]$. 
Then $(R,\a,\b)$ is a BiHom-Lie conformal algebra.
\medskip

  \noindent{\bf Example 2.3.} Let $(L, [\c,\c], \a, \b)$ be a regular BiHom-Lie algebra.
   Denote by $\hat{L}=L\o \mathbb{C}[t, t^{-1}]$ the affization of $L$ with
  \begin{eqnarray*}
 [u\o t^m, v\o t^n]=[u,v]\o t^{m+n},
  \end{eqnarray*}
  for any $u,v\in L, m,n\in \mathbb{Z}$.
  Extend $\a,\b$ to $\hat{L}$ by $\a(u\o t^{m})=\a(u)\o t^{m}$ and $\b(u\o t^{m})=\b(u)\o t^{m}$. 
  Then $(\hat{L}, [\c,\c], \a, \b)$ is a BiHom-Lie algebra. By simple verification, we get a BiHom-Lie conformal algebra $R=\mathbb{C}[\partial]L$ with
\begin{eqnarray*}
[u_{\lambda}v]=[u,v], ~~
\a(f(\partial)u)=f(\partial)\a(u), \b(f(\partial)u)=f(\partial)\b(u), 
\end{eqnarray*}
 for any $u,v\in L.$
\medskip

 \noindent{\bf Definition 2.4.} A module $(M, \a_M, \b_M)$ over over a BiHom-Lie conformal algebra $(R, \a,\b)$ is a $\mathbb{C}[\partial]$-
module endowed with a $\mathbb{C}$-linear map $\a_M, b_M$ and a $\mathbb{C}$-bilinear map $R\o M\rightarrow M[\lambda], a\o v\mapsto a_\lambda v$, 
such that for $a,b\in R, v\in M$:
 \begin{eqnarray*}
\a_M \circ \b_M=\b_M\circ \a_M,\\
\a_M(a\c m)=\a(a)\c \a_M(m),\\
\b_M(a\c m)=\b(a)\c \b_M(m),\\
\a\b(a)_\lambda(b_\mu v)-\b(b)_{\mu}(\a(a)_\lambda v)=[\b(a)_\lambda b]_{\lambda+\mu}\b_{M}(v),\\
(\partial a)_{\lambda}v=-\lambda(a_\lambda v),a_\lambda(\partial v)=(\partial+\lambda)a_{\lambda} v,\\
\b_M \circ \partial =\partial \circ \b_{M}, \a_M \circ \partial =\partial \circ \a_{M},\\
\a_M(a_\lambda v)=\a(a)_\lambda(\a_M(v)),\b_M(a_\lambda v)=\b(a)_\lambda(\b_M(v)),
 \end{eqnarray*}

\noindent{\bf Example  2.5.} Let $(R, \a,\b)$ be a BiHom-Lie conformal algebra. Then $(R, \a,\b)$ is an $R$-module
under the adjoint diagonal action, namely, $a_\lambda b := [a_\lambda b],$  for any $a,b\in R$.

 \noindent{\bf Proposition 2.6.} Let $(R, \a,\b)$ be a regular BiHom-Lie conformal algebra and $(M,\a_M,\b_M)$
an $R$-module. Assume that $\a$ and $\b_M$ are bijective. Define a $\lambda$-bracket $[\c_{\lambda}\c]$ on $R\oplus M$ by
\begin{eqnarray*}
[(a+u)_\lambda (b+v)]_{M}=[a_\lambda b] +a_\lambda v-\a^{-1}\b(b)_{-\partial-\lambda}\a_M\b_M^{-1}(u),
\end{eqnarray*}
for any $a,b\in R$ and $u,v\in M.$
Define $\a+\a_M, \b+\b_M: R\oplus M\rightarrow R\oplus M$ by
 \begin{eqnarray*}
 (\a+\a_M)(a+u)=\a(a)+\a_M(u),~~(\b+\b_M)(a+u)=\b(a)+\b_M(u).
  \end{eqnarray*}
  Then $(R\oplus M, \a+\a_M, \b+\b_M)$ is a regular  BiHom-Lie conformal algebra.
\medskip

{\bf Proof.} Note that $R\oplus M$ is equipped with a $\mathbb{C}[\partial]$-module structure via
\begin{eqnarray*}
\partial(a+u)=\partial(a)+\partial(u),~~\forall a\in R, u\in M.
\end{eqnarray*}
With this, it is easy to see that $(\a+\a_M)\circ \partial =\partial\circ (\a+\a+M),(\b+\b_M)\circ \partial =\partial\circ (\b+\b_M) $ and $(\a+\a_M)([(a+u)_\lambda(b+v)]_M)=[((\a+\a_M)(a+u))_\lambda(\a+\a_M)(b+v)]_M,  (\b+\b_M)([(a+u)_\lambda(b+v)]_M)=[((\b+\b_M)(a+u))_\lambda(\b+\b_M)(b+v)]_M$, 
for any $a,b\in R$ and $u,v\in M.$ 
For (2.1), we have
\begin{eqnarray*}
&&[\partial(a+u)_\lambda (b+v)]_{M}
=[(\partial a+\partial u)_\lambda (b+v)]_{M}\\
&=& [\partial a_{\lambda}b]+(\partial a)_\lambda v-\a^{-1}\b(b)_{-\partial-\lambda}\partial \a_M\b_M^{-1}(u)\\
&=& -\lambda [a_{\lambda}b]-\lambda a_\lambda v-(\partial-\lambda-\partial)\a^{-1}\b(b)_{-\partial-\lambda} \a_M\b_M^{-1}(u)\\
&=& -\lambda([a_{\lambda}b]+ a_\lambda v-\a^{-1}\b(b)_{-\partial-\lambda} \a_M\b_M^{-1}(u))\\
&=& -\lambda[(a+u)_\lambda (b+v)]_{M},\\
&&[(a+u)_\lambda \partial(b+v)]_{M}
= [(a+u)_\lambda (\partial b+ \partial v)]_{M}\\
&=& [a_\lambda \partial b]+a_\lambda \partial v-\partial \a^{-1}\b(b)_{-\partial-\lambda} \a_M\b_M^{-1}(u)\\
&=& (\partial+\lambda)[[a_{\lambda}b]]+(\partial+\lambda)a_\lambda v-(\partial+\lambda)\a^{-1}\b(b)_{-\partial-\lambda} \a_M\b_M^{-1}(u)\\
&=&(\partial+\lambda)([a_{\lambda}b]+ a_\lambda v-\a^{-1}\b(b)_{-\partial-\lambda} \a_M\b_M^{-1}(u))\\
&=&(\partial+\lambda) [(a+u)_\lambda (b+v)]_{M}.
\end{eqnarray*}
as desired. For (2.2), we calculate
\begin{eqnarray*}
&&[(\b(b)+\b_{M}(v))_{-\partial-\lambda}(\a(a)+\a_M(u))]_M\\
&=&[\b(b)_{-\partial-\lambda}\a(a)]+\b(b)_{-\partial-\lambda}\a_M(u)-\b(a)_\lambda \a_M(v)\\
&=& -[\b(a)_\lambda\a(b)]-\b(a)_\lambda \a_M(v)+\b(b)_{-\partial-\lambda}\a_M(u)\\
&=& -[(\b(a)+\b_M(u))_\lambda(\a(b)+\a_M(v))]_M.
\end{eqnarray*}

For (2.3), we have
\begin{eqnarray*}
&&[(\a+\a_M)(\b+\b_M)(a+u)_{\lambda}[(b+v)_{\mu}(c+w)]_M]_M\nonumber\\
&=&  [(\a\b(a)+\a_M\b_M(u))_{\lambda}[(b+v)_{\mu}(c+w)]_M]_M\nonumber
\end{eqnarray*}
\begin{eqnarray*}
&=& [(\a\b(a)+\a_M\b_M(u))_{\lambda}([b_\mu c]+b_\mu w-\a^{-1}\b(c)_{-\partial-\mu}\a_M\b^{-1}_M(v))]_M\nonumber\\
&=&[\a\b(a)_\lambda [b_\mu c]]+\a\b(a)_\lambda(b_\mu w)-\a\b(a)_\lambda(\a^{-1}\b(c)_{-\partial-\mu}\a_M\b^{-1}_M(v))\nonumber\\
&&-\a^{-1}\b([b_\mu c])_{-\partial-\lambda}\a^2_M(u),\\
&&[(\b+\b_M)(b+v)_\mu[(\a+\a_M)(a+u)_{\lambda}(c+w)]_{M}]_M\nonumber\\
&=&[\b(b)_\mu [\a(a)_\lambda c]]+\b(b)_\mu(\a(a)_\lambda w)-\b(b)_\mu(\a^{-1}\b(c)_{-\partial-\mu}\a_M^2\b^{-1}_M(u))\nonumber\\
&&-\a^{-1}\b([\a(a)_\lambda c])_{-\partial-\mu}\a_M(v),\\
&&[[(\b+\b_M)(a+u)_\lambda (b+v)]_{M\lambda+\mu}(\b+\b_M)(c+w)]_M\nonumber\\
&=&[[\b(a)_\lambda b ]_{\lambda+\mu} \b(c)]+[\b(a)_\lambda b]_{\lambda+\mu}\b_M(w)-\a^{-1}\b^2(c)_{-\partial-\lambda-\mu}\a_M\b_M^{-1}(\b(a)_\lambda v)\nonumber\\
&&-\a^{-1}\b^2(c)_{-\partial-\lambda-\mu}\a_M\b^{-1}_M(\a^{-1}\b(b)_{-\partial-\lambda}\a_{M}(u)).\nonumber
\end{eqnarray*}
Since $(M,\a_M,\b_M)$ is an $R$-module, we have
\begin{eqnarray*}
\a\b(a)_\lambda(b_\mu v)-\b(b)_{\mu}(\a(a)_\lambda v)=[\b(a)_\lambda b]_{\lambda+\mu}\b_{M}(v),
\end{eqnarray*}
This implies
\begin{eqnarray*}
&&[(\a+\a_M)(\b+\b_M)(a+u)_{\lambda}[(b+v)_{\mu}(c+w)]_M]_M\\
&=& [(\b+\b_M)(b+v)_\mu[(\a+\a_M)(a+u)_{\lambda}(c+w)]_{M}]_M\\
&&+[[(\b+\b_M)(a+u)_\lambda (b+v)]_{M\lambda+\mu}(\b+\b_M)(c+w)]_M.
\end{eqnarray*}
Therefore,  $(R\oplus M, \a+\a_M, \b+\b_M)$ is a regular BiHom-Lie conformal algebra. \hfill $\square$
\medskip

In the following we aim to develop the cohomology theory of regular BiHom-Lie conformal algebras.
To do this, we need the following concept.
\medskip

 \noindent{\bf Definition 2.7.}
An $n$-cochain ($n\in \mathbb{Z}_{\geq0}$) of a regular BiHom-Lie conformal algebra $R$ with coefficients
 in a module $(M,\alpha_{M},\b_{M} )$ is a $\mathbb{C}$-linear map
\begin{eqnarray*}
\g:R^n\rightarrow M[\lambda_1,... ,\lambda_n],~~
(a_1,...,a_n)\mapsto \g_{\lambda_1,...,\lambda_n}(a_1,... , a_n),
\end{eqnarray*}
where $M[\lambda_1,...,\lambda_n]$ denotes the space of polynomials with coefficients in $M$, satisfying the following conditions:

(1) Conformal antilinearity:
\begin{eqnarray*}
\g_{\lambda_1,...,\lambda_n}(a_1,... ,\partial a_i,... ,a_n)=-\lambda_i\g_{\lambda_1,...,\lambda_n}(a_1,... ,a_i,... ,a_n).
\end{eqnarray*}

(2) Skew-symmetry:
\begin{eqnarray*}
&&\gamma(a_{1},..., \b(a_{i+1}), \a(a_{i}),..., a_{n})=-\gamma(a_{1}, ..., \b(a_{i}), \a(a_{i+1}), ... , a_{n}).
\end{eqnarray*}

(3) Commutativity:
\begin{eqnarray*}
\g\circ \a=\a_M\circ \g, \g\circ \b=\b_M\circ \g,
\end{eqnarray*}

Let $R^{\otimes 0} = \mathbb{C}$ as usual so that a $0$-cochain
 is an element of $M$. Define a differential $d$ of a cochain $\g$ by
\begin{eqnarray*}
&&(d\gamma)_{\lambda_1,...,\lambda_{n+1}}(a_{1},...,a_{n+1})\nonumber\\
&=&\sum_{i=1}^{ n+1}(-1)^{i+1}\alpha\b^{n-1}(a_{i})_{\lambda_i}\gamma_{\lambda_1,...,\hat{\lambda}_{i},...,\lambda_{n+1}}(a_{1},..., \hat{a}_{i},..., a_{n+1})+\\
  &&\sum_{1\leq i<j\leq n+1} (-1)^{i+j}\gamma_{\lambda_i+\lambda_j, \lambda_1,...,\hat{\lambda}_i,...\hat{\lambda}_j,...,\lambda_{n+1}}([\a^{-1}\b(a_{i})_{\lambda_i}a_{j}], \b(a_{1}),... ,\hat{a}_{i},... , \hat{a}_{j},...,\b(a_{n+1})),
\end{eqnarray*}
where $\g$ is extended linearly over the polynomials in
$\lambda_i$. In particular, if $\g$ is a 0-cochain, then $(d\g)_\lambda a=a_\lambda \g$.

 \noindent{\bf Proposition 2.8.} $d\g$ is a cochain and $d^2=0$.

 {\bf Proof.} Let $\g$ be an $n$-cochain.  It is easy to check that $d$ satisfies conformal antilinearity, skew-symmetry  and commutativity.
  That is, $d$ is an $(n+1)$-cochain.
 Next we will check that $d^2=0$.  In fact, we have
 \begin{eqnarray}
&& (d^2\g) _{\lambda_1,...,\lambda_{n+2}}(a_{1},...,a_{n+2})\nonumber\\
&=& \sum_{i=1}^{ n+1}(-1)^{i+1}\alpha\b^{n}(a_{i})_{\lambda_i}(d\g)_{\lambda_1,...,\hat{\lambda}_{i},...,\lambda_{n+2}}(a_{1},..., \hat{a}_{i},..., a_{n+2})\nonumber\\
  &&+\sum_{1\leq i<j\leq n+1} (-1)^{i+j}(d\g)_{\lambda_i+\lambda_j, \lambda_1,...,\hat{\lambda}_i,...\hat{\lambda}_j,...,\lambda_{n+2}}([\a^{-1}\b(a_{i})_{\lambda_i}a_{j}], \b(a_{1}),... ,\hat{a}_{i},... , \hat{a}_{j},...,\b(a_{n+2}))\nonumber\\
  &=&\sum_{i=1}^{n+2}\sum_{j=1}^{i-1}(-1)^{i+j}\alpha\b^{n}(a_{i})_{\lambda_i}(\alpha\b^{n-1}(a_{j})_{\lambda_j}\g_{\lambda_1,...,\hat{\lambda}_{j,i},...,\lambda_{n+2}}(a_{1},..., \hat{a}_{j,i},..., a_{n+2})\\
  &&\sum_{i=1}^{n+2}\sum_{j=i+1}^{n+2}(-1)^{i+j+1}\alpha\b^{n}(a_{i})_{\lambda_i}(\alpha\b^{n-1}(a_{j})_{\lambda_j}\g_{\lambda_1,...,\hat{\lambda}_{i,j},...,\lambda_{n+2}}(a_{1},..., \hat{a}_{i,j},..., a_{n+2})\\
  &&+\sum_{i=1}^{n+2}\sum_{1\leq j< k<i}^{n+2}(-1)^{i+j+k+1}\alpha\b^{n}(a_{i})_{\lambda_i}\g_{\lambda_j+\lambda_k, \lambda_1,...,\hat{\lambda}_{j,k,i},...,\lambda_{n+2}}\nonumber\\
  &&([\a^{-1}\b(a_{j})_{\lambda_j}a_{k}], \b(a_{1}),... ,\hat{a}_{j,k,i},...,\b(a_{n+2}))\\
  &&+\sum_{i=1}^{n+2}\sum_{1\leq j< i<k}^{n+2}(-1)^{i+j+k}\alpha\b^{n}(a_{i})_{\lambda_i}\g_{\lambda_j+\lambda_k, \lambda_1,...,\hat{\lambda}_{j,i,k},...,\lambda_{n+2}}\nonumber\\
  &&([\a^{-1}\b(a_{j})_{\lambda_j}a_{k}], \b(a_{1}),... ,\hat{a}_{j,i,k},...,\b(a_{n+2}))\\
   &&+\sum_{i=1}^{n+2}\sum_{1\leq i< j<k}^{n+2}(-1)^{i+j+k+1}\alpha\b^{n}(a_{i})_{\lambda_i}\g_{\lambda_j+\lambda_k, \lambda_1,...,\hat{\lambda}_{i,j,k},...,\lambda_{n+2}}\nonumber\\
  &&([\a^{-1}\b(a_{j})_{\lambda_j}a_{k}], \b(a_{1}),... ,\hat{a}_{i,j,k},...,\b(a_{n+2}))
     \end{eqnarray}
   \begin{eqnarray}
  &&+\sum_{1\leq i< j}^{n+2}\sum_{k=1}^{i-1}(-1)^{i+j+k}\alpha\b^{n}(a_{k})_{\lambda_k}\g_{\lambda_i+\lambda_j, \lambda_1,...,\hat{\lambda}_{k, i,j},...,\lambda_{n+2}}\nonumber\\
  &&([\a^{-1}\b(a_{i})_{\lambda_i}a_{j}], \b(a_{1}),... ,\hat{a}_{k, i,j},...,\b(a_{n+2}))\\
  &&+\sum_{1\leq i< j}^{n+2}\sum_{k=i+1}^{j-1}(-1)^{i+j+k+1}\alpha\b^{n}(a_{k})_{\lambda_k}\g_{\lambda_i+\lambda_j, \lambda_1,...,\hat{\lambda}_{i, k,j},...,\lambda_{n+2}}\nonumber\\
  &&([\a^{-1}\b(a_{i})_{\lambda_j}a_{j}], \b(a_{1}),... ,\hat{a}_{i, k,j},...,\b(a_{n+2}))\\
  &&+\sum_{1\leq i< j}^{n+2}\sum_{k=j+1}^{n+2}(-1)^{i+j+k}\alpha\b^{n}(a_{k})_{\lambda_k}\g_{\lambda_i+\lambda_j, \lambda_1,...,\hat{\lambda}_{i, j,k},...,\lambda_{n+2}}\nonumber\\
  &&([\a^{-1}\b(a_{i})_{\lambda_i}a_{j}], \b(a_{1}),... ,\hat{a}_{i,j,k},...,\b(a_{n+2}))\\
  &&+\sum_{1\leq i< j}^{n+2}(-1)^{i+j}\alpha\b^{n-1}([\a^{-1}\b(a_{i})_{\lambda_i}a_{j})_{\lambda_i+\lambda_j}\g_{\lambda_1,...,\hat{\lambda}_{j},...\hat{\lambda}_{i},...,
  \lambda_{n+2}}(\b(a_{1}),... ,\hat{a}_{j},...,\hat{a}_{i},...,\b(a_{n+2})~~~~~~~~~\\
  &&+\sum^{n+2}_{distinct i,j,k,l,i< j,k<l}(-1)^{i+j+k+l}sign\{i,j, k, l\}\nonumber\\
  &&\g_{\lambda_k+\lambda_l, \lambda_i+\lambda_j,\lambda_1,...,\hat{\lambda}_{i,j,k,l},...,\lambda_{n+2}}(\b[\a^{-1}\b(a_{k})_{\lambda_k}a_{l}],
  \b[\a^{-1}\b(a_{i})_{\lambda_i}a_{j}],...,\hat{a}_{i,j,k,l},..., \b^2(a_{n+2}))~~~~~~~~~\\
   &&\sum^{n+2}_{i,j,k=1,i< j,k\neq i,j}(-1)^{i+j+k+l}sign\{i,j, k\}\nonumber\\
  &&\g_{\lambda_i+\lambda_k+\lambda_j,\lambda_1,...,\hat{\lambda}_{i,j,k},...,\lambda_{n+2}}([[\a^{-1}\b(a_i)_{\lambda_i}a_{j}]_{\lambda_i+\lambda_j},\b(\a_k)],
  \b^2(a_1),...,\hat{a}_{i,j,k},..., \b^2(a_{n+2})),~~~~~~~~~
 \end{eqnarray}
where $sign\{i_1, ... , i_p\}$ is the sign of the permutation putting the indices in increasing
order and $\hat{a}_{i,j}$, $a_i, a_j,... $ are omitted.
\medskip

It is obvious that (2.6) and (2.11) summations cancel each other. The same is true
for (2.7) and (2.10), (2.8) and (2.9). The BiHom-Jacobi identity implies (2.14) = 0 and
the skew-symmetry implies  (2.13)=0. Because $(M, \a_M,\b_M)$ is an $R$-module, it follows that
\begin{eqnarray*}
\a\b(a)_\lambda(b_\mu v)-\b(b)_{\mu}(\a(a)_\lambda v)=[\b(a)_\lambda b]_{\lambda+\mu}\b_{M}(v).
\end{eqnarray*}
Since $\g\circ \a=\a_M\circ \g, \g\circ \b=\b_M\circ \g$, we have (2.4), (2.5) and (2.12) summations cancel. 
So $d^2\g=0$ and the proof is completed. \hfill $\square$
\medskip

 Thus the cochains of a BiHom-Lie conformal algebra $R$ with coefficients in a module $M$ form a
complex, which is denoted by
 \begin{eqnarray*}
\widetilde{C}^{\bullet}_{\a,\b} =\widetilde{C}^{\bullet}_{\a,\b}(R,M)=\bigoplus_{n\in \mathbb{Z}_{+}}C^{n}_{\a,\b}(R,M).
 \end{eqnarray*}
This complex is called the basic complex for the $R$-module $(M,\a_M,\b_M)$.

 \section{ Deformations  and  Nijenhuis operators of  BiHom-Lie\\ conformal algebras}
\def\theequation{\arabic{section}. \arabic{equation}}
\setcounter{equation} {0}
 In this section, we give the definition of a deformation and show that the deformation generated
by a 2-cocycle Nijenhuis operator is trivial.
\medskip

 Let $R$ be a regular BiHom-Lie conformal algebra. Define
$
a_\lambda b=[\a^{s}(a)_\lambda b],$ for any $a,b\in R$.
 Set $\g\in \widetilde{C}^n_{\a,\b}(R,R_s)$.  Define an operator $d_s: \widetilde{C}^{n}_{\a,\b}(R,R_s)\rightarrow \widetilde{C}^{n+1}_{\a,\b}(R,R_s)$ by
\begin{eqnarray*}
&&(d_s\gamma)_{\lambda_1,...,\lambda_{n+1}}(a_{1},..., a_{n+1})\nonumber\\
&=&\sum_{i=1}^{ n+1}(-1)^{i+1}[\a^{s+1}\b^{n-1}(a_{i})_{\lambda_i}\gamma_{\lambda_1,...,\hat{\lambda}_{i},...,\lambda_{n+1}}(a_{1},..., \hat{a}_{i},..., a_{n+1})]+\nonumber\\
  &&\sum_{1\leq i<j\leq n+1}(-1)^{i+j} \gamma_{\lambda_i+\lambda_j, \lambda_1,...,\hat{\lambda}_i,...\hat{\lambda}_j,...,\lambda_{n+1}}([\a^{-1}\b(a_{i})_{\lambda_i}a_{j}], \b(a_{1}),... ,\hat{a}_{i},... , \hat{a}_{j},...,\b(a_{n+1})).
\end{eqnarray*}
 Obviously, the operator $d_s$ is induced from the differential $d$. Thus $d_s$ preserves the space of cochains
and satisfies $d_s^2=0$. In the following the complex $\widetilde{C}^{\bullet}_{\a,\b}(R, R_s)$ is assumed to be associated with the
differential $d_s$.
\medskip

Taking $s=-1$, let $\psi\in C^{2}(R,R)_{\overline{0}}$ be a bilinear operator commuting with $\a$ and $\b$,  we consider a $t$-parameterized family of bilinear operations on $R$,
\begin{eqnarray}
[a_\lambda b]_t=[a_\lambda b]+t\psi_{\lambda, -\partial-\lambda}(a,b),~~~\forall a,b\in R.
\end{eqnarray}
If $[\c_{\lambda}\c]$ endows $(R,[\c_{\lambda}\c], \a,\b)$ with a BiHom-Lie conformal algebra structure, we say that $\psi$  generates
a deformation of the BiHom-Lie conformal algebra $R$. It is easy to see that $[\c_{\lambda}\c]$ satisfies (2.1).
 The skew-symmetry of $[\c_{\lambda}\c]$ means that
 \begin{eqnarray*}
&&[\b(a)_\lambda \a(b)]_t=[\b(a)_\lambda \a(b)]+t\psi_{\lambda, -\partial-\lambda}(\b(a),\a(b)),\\
&&[\b(b)_{-\partial-\lambda} \a(a)]_t=[\b(b)_{-\partial-\lambda}\a(a)]+t\psi_{-\partial-\lambda, \lambda}(\b(b),\a(a)).
 \end{eqnarray*}
Then $[\b(a)_\lambda \a(b)]_t=-[\b(b)_{-\partial-\lambda} \a(a)]_t$   if and only if
 \begin{eqnarray}
\psi_{\lambda, -\partial-\lambda}(\b(a),\a(b))=-\psi_{-\partial-\lambda, \lambda}(\b(b),\a(a)).
 \end{eqnarray}
 If it is true for (2.3), by expanding the BiHom-Jacobi identity for $[\c_{\lambda}\c]$, we have
 \begin{eqnarray*}
 &&[\a\b(a)_\lambda[b_{\mu}c]]+t([\a\b(a)_\lambda(\psi_{\mu,-\partial-\mu}(b,c)]+\psi_{\lambda,-\partial-\lambda}(\a\b(a),[b_{\mu}c]))\\
 &&+t^2\psi_{\lambda,-\partial-\lambda}(\a\b(a),\psi_{\mu,-\partial-\mu}(b,c))\\
 &=& [\b(b)_\mu[\a(a)_{\lambda}c]]+t([\b(b)_\mu(\psi_{\lambda,-\partial-\lambda}(\a(a),c))]+\psi_{\mu,-\partial-\mu}(\b(b),[\a(a)_{\lambda}c]))\\
 &&+t^2\psi_{\mu,-\partial-\mu}(\b(b),\psi_{\lambda,-\partial-\lambda}(\a(a),c))+[[\b(a)_{\lambda}b]_{\lambda+\mu}\b(c)]\\
 &&+t([(\psi_{\lambda,-\partial-\lambda}(\b(a),b))_{\lambda+\mu}\b(c)]+ \psi_{\lambda+\mu,-\partial-\lambda-\mu}([\b(a)_\lambda b],\b(c)))\\
 &&+t^2\psi_{\lambda+\mu,-\partial-\lambda-\mu}(\psi_{\lambda,-\partial-\lambda}(\b(a),b), \b(c)).
 \end{eqnarray*}
This is equivalent to the following conditions
\begin{eqnarray}
 &&\psi_{\lambda,-\partial-\lambda}(\a\b(a),\psi_{\mu,-\partial-\mu}(b,c))\nonumber\\
 &=&\psi_{\mu,-\partial-\mu}(\b(b),\psi_{\lambda,-\partial-\lambda}(\a(a),c))
  + \psi_{\lambda+\mu,-\partial-\lambda-\mu}(\psi_{\lambda,-\partial-\lambda}(\b(a),b),\b(c)),\\
&&[\a\b(a)_\lambda(\psi_{\mu,-\partial-\mu}(b,c)]+\psi_{\lambda,-\partial-\lambda}(\a\b(a),[b_{\mu}c])\nonumber\\
&=& [\b(b)_\mu(\psi_{\lambda,-\partial-\lambda}(\a(a),c))]+\psi_{\mu,-\partial-\mu}(\b(b),[\a(a)_{\lambda}c])\nonumber\\
&&+[(\psi_{\lambda,-\partial-\lambda}(\b(a),b))_{\lambda+\mu}\b(c)]+ \psi_{\lambda+\mu,-\partial-\lambda-\mu}([\b(a)_\lambda b],\b(c)).
\end{eqnarray}
Obviously, (3.3) and (3.4) mean  that $\psi$  define a BiHom-Lie conformal algebra structure on $R$, as required.
\hfill $\square$
\medskip

A deformation is said to be trivial if there is a linear operator $f\in \widetilde{C}_{\a,\b}^{1}(R,R_{-1})$ such that
\begin{eqnarray}
T_t([a_\lambda b]_t)=[T_t(a)_\lambda T_t(b)],~~~
\end{eqnarray}
 for any $a,b\in R$, where $T_t=id+tf$.
\medskip

 \noindent{\bf Definition  3.1.} A linear operator $f\in \widetilde{C}_{\a,\b}^{1}(R,R_{-1})$ is a Nijienhuis operator
if  \begin{eqnarray}
    [f(a)_\lambda f(b)]=f([a_\lambda b]_{N}),
    \end{eqnarray}
 for any $a,b\in R$, where the bracket $[\c,\c]_{N}$ is defined by
\begin{eqnarray}
[a_\lambda b]_{N}=[f(a)_\lambda b]+[a_\lambda f(b)]-f([a_\lambda b]).
\end{eqnarray}

\noindent{\bf Theorem  3.2.} Let $(R,\a,\b)$ be a regular BiHom-Lie conformal algebra and $f\in \widetilde{C}_{\a,\b}^{1}(R,R_{-1})$
a Nijienhuis operator. Then a deformation of $(R,\a,\b)$) can be obtained by putting
\begin{eqnarray}
\psi_{\lambda, -\partial-\lambda}(a,b)=[a_\lambda b]_{N},
\end{eqnarray}
for any $a,b\in R$. Furthermore, this deformation is trivial.

{\bf Proof.}  To show that $\psi$  generates a deformation, we need to verify (3.2), (3.3)
and (3.4). We first prove that $\psi(\b(a),\a(b))=-\psi(\b(b),\a(a)).$ In fact, we have
\begin{eqnarray*}
\psi(\b(a),\a(b))
&=&[a_\lambda b]_N\\
&=&[f(\b(a))_\lambda \a(b)]+[\b(a)_\lambda f(\a(b))]-f[(\b(a)_\lambda\a(b)]\\
&=& -([f(\b(b))_\lambda \a(a)]+[\b(b)_\lambda f(\a(a))]-f[(\b(b)_\lambda\a(a)])\\
&=& -[\b(b),\a(a)]_N\\
&=&-\psi(\b(b),\a(a)),
\end{eqnarray*}
as desired. By (3.5),(3.6) and (3.7), we have
\begin{eqnarray*}
&&\psi_{\lambda,-\partial-\lambda}(\a\b(a),\psi_{\mu,-\partial-\mu}(b,c))+\psi_{\mu,-\partial-\mu}(\b(b),\psi_{-\partial-\lambda,\lambda}(\a^{-1}\b(c),\a^2\b^{-1}(a)))\\
&&+\psi_{-\partial-\lambda-\mu, \lambda+\mu}(\a^{-1}\b^2(c), \psi_{\lambda,-\partial-\lambda}(\a(a),\a\b^{-1}(b)))\\
&=&[f(\a\b(a))_{\lambda}[f(b)_{\mu}c]]+[f(\a\b(a))_{\lambda}[b_\mu f(c)]]-[f(\a\b(a))_{\lambda}f([b_{\mu}c])] \\
&&+[\a\b(a)_{\lambda}[f(b)_\mu f(c)]]-f([\a\b(a)_{\lambda}[f(b)_\mu c]])-f([\a\b(a)_{\lambda}[b_\mu f(c)]])\\
&&+f([\a\b(a)_{\lambda}f([b_\mu c])])+[f(\b(b))_{\mu}[f(\a^{-1}\b(c))_{-\partial-\lambda}\a^2\b^{-1}(a)]]\\
&&+[f(\b(b))_{\mu}[\a^{-1}\b(c)_{-\partial-\lambda} f(\a^2\b^{-1}(a))]]-[f(\b(b))_{\mu}f([\a^{-1}\b(c)_{-\partial-\lambda}\a^2\b^{-1}(a)])]\\
&&+[\b(b)_{\mu}[f(\a^{-1}\b(c))_{-\partial-\lambda} f(\a^2\b^{-1}(a))]]-f([\b(b)_{\mu}[f(\a^{-1}\b(c))_{-\partial-\lambda} \a^2\b^{-1}(a)]])\\
&&-f([\b(b)_{\mu}[\a^{-1}\b(c)_{-\partial-\lambda} f(\a^2\b^{-1}(a))]])+[f(\b(b))_{\mu}f([\a^{-1}\b(c)_{-\partial-\lambda} \a^2\b^{-1}(a)])])\\
&&+[f(\a^{-1}\b^2(c))_{-\partial-\lambda-\mu}[f(\a(a))_{\lambda}\a\b^{-1}(b)]]+[f(\a^{-1}\b^2(c))_{-\partial-\lambda-\mu}[\a(a)_\lambda f(\a\b^{-1}(b))]]\\
&&-[f(\a^{-1}\b^2(c))_{-\partial-\lambda-\mu}f([\a(a)_{\lambda}\a\b^{-1}(b)])]+[\a^{-1}\b^2(c)_{-\partial-\lambda-\mu}[f(\a(a))_{\lambda}f(\a\b^{-1}(b))]]\\
&&-f([\a^{-1}\b^2(c)_{-\partial-\lambda-\mu}[f(\a(a))_{\lambda}\a\b^{-1}(b)]])-f([\a^{-1}\b^2(c)_{-\partial-\lambda-\mu}[\a(a)_{\lambda}f(\a\b^{-1}(b))]])\\
&&+f([\a^{-1}\b^2(c)_{-\partial-\lambda-\mu}f([\a(a)_{\lambda}\a\b^{-1}(b)])]).
\end{eqnarray*}
Since $f$ is a Nijenhuis operator, it follows that 
\begin{eqnarray*}
&&-[f(\a\b(a))_{\lambda}f([b_{\mu}c])]+f([\a\b(a)_{\lambda}f([b_\mu c])])\\
&=&-f([f(\a\b(a))_{\lambda}[b_{\mu}c]])+f^2([\a\b(a)_{\lambda}[b_\mu c]]),\\
&& -[f(\b(b))_{\mu}f([\a^{-1}\b(c)_{-\partial-\lambda}\a^2\b^{-1}(a)])]+f([\b(b)_{\mu}f([\a^{-1}\b(c)_{-\partial-\lambda} \a^2\b^{-1}(a)])])\\
&=&-f([f(\b(b))_{\mu}[\a^{-1}\b(c)_{-\partial-\lambda}\a^2\b^{-1}(a)]])+f^2([\b(b)_{\mu}[\a^{-1}\b(c)_{-\partial-\lambda} \a^2\b^{-1}(a)]]),\\
&& -[f(\a^{-1}\b^2(c))_{-\partial-\lambda-\mu}f([\a(a)_{\lambda}\a\b^{-1}(b)])]+f([\a^{-1}\b^2(c)_{-\partial-\lambda-\mu}f([\a(a)_{\lambda}\a\b^{-1}(b)])])\\
&=&-f[f(\a^{-1}\b^2(c))_{-\partial-\lambda-\mu}[\a(a)_{\lambda}\a\b^{-1}(b)]]+f^2([\a^{-1}\b^2(c)_{-\partial-\lambda-\mu}[\a(a)_{\lambda}\a\b^{-1}(b)]]).
\end{eqnarray*}
Note that
\begin{eqnarray*}
[\a\b(a)_\lambda[b_\mu c]]+[\a^{-1}\b^2(c)_{-\partial-\lambda-\mu}[\a(a)_{\lambda}\a\b^{-1}(b)]]+[\b(b)_{\mu}[\a^{-1}\b(c)_{-\partial-\lambda}\a^2\b^{-1}(a)]]=0.
\end{eqnarray*}
Thus
\begin{eqnarray*}
&&\psi_{\lambda,-\partial-\lambda}(\a\b(a),\psi_{\mu,-\partial-\mu}(b,c))+\psi_{\mu,-\partial-\mu}(\b(b),\psi_{-\partial-\lambda,\lambda}(\a^{-1}\b(c),\a^2\b^{-1}(a)))\\
&&+\psi_{-\partial-\lambda-\mu, \lambda+\mu}(\a^{-1}\b^2(c), \psi_{\lambda,-\partial-\lambda}(\a(a),\a\b^{-1}(b)))=0.
\end{eqnarray*}
In the similar way, we can check that (3.4). This proves that $\psi$ generates a deformation of the regular BiHom-Lie conformal algebra $(R, \a,\b)$.

Let $T_t=id+tf$. By (3.8), we have
\begin{eqnarray}
T_t([a_\lambda b]_t)&=&(id+tf)([a_\lambda b]+t\psi_{\lambda,-\partial-\lambda}(a,b))\nonumber\\
&=& (id+tf)([a_\lambda b]+t[a_\lambda b]_N)\nonumber\\
&=& [a_\lambda b]+t([a_\lambda b]_N+f([a_\lambda b]))+t^2f([a_\lambda b]_N).
\end{eqnarray}
On the other hand, we have
\begin{eqnarray}
[T_t(a)_\lambda T_t(b)]&=&[(a+tf(a))_\lambda(b+tf(b))]\nonumber\\
&=& [a_\lambda b]+t([f(a)_\lambda b]+[a_\lambda f(b)])+t^2[f(a)_\lambda f(b)].
\end{eqnarray}
Combining (3.9) with (3.10), it follows that $T_t([a_\lambda b]_t)=[T_t(a)_\lambda T_t(b)]$. Therefore the deformation
is trivial.\hfill $\square$

\section{Derivations of multiplicative  BiHom-Lie conformal algebras }
\def\theequation{\arabic{section}. \arabic{equation}}
\setcounter{equation} {0}

For convenience, we denote by $\mathcal{A}$ the ring $\mathbb{C}[\partial]$ of polynomials in the indeterminate $\partial$.

\noindent{\bf Definition 4.1.}$^{\cite{Kac98}}$ A conformal linear map between $\mathcal{A}$-modules $V$ and $W$ is a linear map
$\phi: V\rightarrow \mathcal{A}[\lambda]\o_{\mathcal{A}}W$ such that
\begin{eqnarray}
\phi(\partial v)=(\partial+\lambda)(\phi v).
\end{eqnarray}

We will often abuse the notation by writing $\phi: V \rightarrow W$, it is clear from the
context that $\phi$ is a conformal linear map. We will also write $\phi_{\lambda}$ instead of $\phi$ to emphasize
the dependence of $\phi$ on $\lambda$.

The set of all conformal linear maps from $V$ to $W$ is denoted by $Chom(V,W)$ and it is
made into an $\mathcal{A}$-module via
\begin{eqnarray}
(\partial\phi)_{\lambda} v=-\lambda
\end{eqnarray}
We will write $Cend(V)$ for $Chom(V, V)$.

\noindent{\bf Definition 4.2.} Let $(R,\a,\b)$ be a multiplicative BiHom-Lie conformal algebra. Then a
conformal linear map $D_\lambda: R\rightarrow R$ is called an $\a^k\b^l$-derivation of $(R,\a,\b)$ if
\begin{eqnarray}
D_\lambda \circ \a=\a\circ D_\lambda, ~~~D_\lambda \circ \b=\b\circ D_\lambda,\nonumber\\
D_\lambda([a_\mu b])=[D_\lambda(a)_{\lambda+\mu}\a^{k}\b^{l}(b)]+[\a^{k}\b^{l}(a)_{\mu}D_\lambda(b)].
\end{eqnarray}
Denote by $Der_{\a^{s},\b^{l}}$ the set of $\a^{s}\b^{l}$-derivations of the multiplicative BiHom-Lie conformal
algebra $(R,\a,\b)$. For any $a\in R$ satisfying $\a(a)=a, \b(a)=a$, define $D_{k,l}:R\rightarrow R$ by
\begin{eqnarray*}
D_{k,l}(a)_{\lambda}(b)=[a_\lambda\a^{k+1}\b^{l-1}(b)],~~~\forall b\in R.
\end{eqnarray*}
Then $D_{k,l}(a)$ is an $\a^{k+1}\b^{l}$-derivation, which is called an inner $\a^{k+1}\b^{l}$-derivation. In fact,
\begin{eqnarray*}
D_{k,l}(a)_{\lambda}(\partial b)&=&[a_\lambda\a^{k+1}\b^{l-1}(\partial b)]
= [a_\lambda\partial\a^{k+1}\b^{l-1}( b)]
= (\partial+\lambda)D_{k,l}(a)_{\lambda}(b),
\end{eqnarray*}
\begin{eqnarray*}
&&D_{k,l}(a)_{\lambda}(\a(b))=[a_\lambda\a^{k+2}\b^{l-1}(b)]
= \a[a_\lambda\a^{k+1}\b^{l-1}( b)]
= \a\circ D_{k,l}(a)_{\lambda}(b),\\
&&D_{k,l}(a)_{\lambda}(\b(b))=[a_\lambda\a^{k+1}\b^{l}(b)]
= \b[a_\lambda\a^{k+1}\b^{l-1}( b)]
= \b\circ D_{k,l}(a)_{\lambda}(b),\\
&&D_{k,l}(a)_{\lambda}([b_{\mu}c])= [a_\lambda \a^{k+1}\b^{l-1}([b_{\mu}c])
=  [\a\b(a)_\lambda [\a^{k+1}\b^{l-1}(b)_{\mu}\a^{k+1}\b^{l-1}(c)]\\
&&~~~~~~~~~~~~~~~~~~~= [\b(a)_\lambda\a^{k+1}\b^{l-1}(b)]_{\lambda+\mu}\a^{k+1}\b^{l}(c)]+[\a^{k+1}\b^{l}(b)_{\mu}[\a(a)_{\lambda}\a^{k+1}\b^{l-1}(c)]]\\
&&~~~~~~~~~~~~~~~~~~~=[a_\lambda\a^{k+1}\b^{l-1}(b)]_{\lambda+\mu}\a^{k+1}\b^{l}(c)]+[\a^{k+1}\b^{l}(b)_{\mu}[a_{\lambda}\a^{k+1}\b^{l-1}(c)]]\\
&&~~~~~~~~~~~~~~~~~~~=[D_{k,l}(a)_{\lambda}(b)_{\lambda+\mu}\a^{k+1}\b^{l}(c)]+[\a^{k+1}\b^{l}(b)_{\mu}(D_{k,l}(a)_{\lambda}(c))].
\end{eqnarray*}
Denote by $Inn_{\a^k\b^l} (R)$ the set of inner $\a^k\b^l$-derivations. 
For any $D_{\lambda}\in Der_{\a^k\b^l}(R)$ and $D'_{\mu-\lambda}\in Der_{\a^s\b^t}(R)$, define their commutator $[D_\lambda D']_{\mu}$ by
\begin{eqnarray}
[D_\lambda D']_{\mu}(a)=D_\lambda(D'_{\mu-\lambda}a)-D'_{\mu-\lambda}(D_\lambda a),~~~\forall a\in R.
\end{eqnarray}

\noindent{\bf Lemma 4.3.} For any $D_{\lambda}\in Der_{\a^k\b^l}(R)$ and $D'_{\mu-\lambda}\in Der_{\a^s\b^t}(R)$, we have
\begin{eqnarray*}
[D_\lambda D'] \in Der_{\a^{k+s}\b^{l+t}}(R)[\lambda].
\end{eqnarray*}
{\bf Proof.} For any $a,b\in R$, we have
\begin{eqnarray*}
&&[D_\lambda D']_{\mu}(\partial a)\\
&=&D_{\lambda}(D'_{\mu-\lambda}\partial a)-D'_{\mu-\lambda}(D_\lambda \partial a)\\
&=& D_\lambda((\partial+\mu-\lambda)D'_{\mu-\lambda}a)+D'_{\mu-\lambda}((\mu+\lambda)D_\lambda a)\\
&=& (\partial+\mu)D_{\lambda}(D'_{\mu-\lambda}a)-(\partial+\mu)D'_{\mu-\lambda}(D_{\lambda}a)\\
&=& (\partial+\mu)[D_\lambda D']_{\mu}(a).
\end{eqnarray*}
Similarly ,we have
\begin{eqnarray*}
&&[D_\lambda D']_{\mu}([a_{\g}b])\\
&=& D_\lambda(D'_{\mu-\lambda}[a_{\g}b])-D'_{\mu-\lambda}(D_\lambda[a_{\g}b])\\
&=&D_\lambda([D'_{\mu-\lambda}(a)_{\mu-\lambda+\gamma}\a^s\b^{t}(b)]+[\a^s\b^{t}(a)_{\g}D'_{\mu-\lambda}(b)])\\
&& -D'_{\mu-\lambda}([D_\lambda(a)_{\lambda+\gamma}\a^k\b^{l}(b)]+[\a^k\b^{l}(a)_{\g}D_{\lambda}(b)])\\
&=& [D_\lambda(D'_{\mu-\lambda}(a))_{\mu+\gamma}\a^{k+s}\b^{l+t}(b)]+[\a^{k}\b^{l}(D'_{\mu-\lambda}(a))_{\mu-\lambda+\g}D_{\lambda}(\a^{s}\b^{t}(b)))\\
&&+[D_\lambda(\a^s\b^t(a))_{\lambda+\gamma}\a^k\b^l(D'_{\mu-\lambda}(b))]+[\a^{k+s}\b^{l+t}(a)_\g(D_\lambda(D'_{\mu-\lambda}(b)))]\\
&&-[(D'_{\mu-\lambda}D_\lambda(a))_{\mu+\gamma}\a^{k+s}\b^{l+t}(b)]-[\a^{s}\b^{t}(D_\lambda(s))_{\lambda+\gamma}(D'_{\mu-\lambda}(\a^k\b^l(b)))]\\
&&-[(D'_{\mu-\lambda}(\a^{k}\b^{l}(a)))_{\mu-\lambda+\gamma}\a^{s}\b^{t}(D_{\lambda}(b))]-[\a^{k+s}\b^{l+t}(a)_{\lambda}(D'_{\mu-\lambda}(D_\lambda(b)))]\\
&=&[([D_\lambda D']_\mu a)_{\mu+\g}\a^{k+s}\b^{l+t}(b)]+[\a^{k+s}\b^{l+t}(a)_{\g}([D_\lambda D']_{\mu}b)].
\end{eqnarray*}
Therefore, $[D_\lambda D'] \in Der_{\a^{k+s}\b^{l+t}}(R)[\lambda]$. \hfill $\square$

Define
\begin{eqnarray}
Der(R)=\bigoplus_{k\geq0, l\geq 0}Der_{\a^k\b^l}(R).
\end{eqnarray}

\noindent{\bf Proposition 4.4.} With notations above, $(Der(R), \a',\b')$ is a BiHom-Lie conformal algebra with respect to (4.4),
 where $\a'(D)=D\circ \a, \b'(D)=D\circ \b$.

{\bf Proof.}  By (4.2), $Der(R)$ is a $\mathbb{C}[\partial]$-module. By (4.1), (4.2) and (4.4), it is easy to check that (2.1) and (2.2)
are satisfied. To check the BiHom-Jacobi identity, we compute 
\begin{eqnarray*}
&&[\a'\b'(D)_{\lambda}[D'_{\mu}D'']]_{\theta}(a)\\
&=&(D\circ \a\b)_{\lambda}([D'_{\mu}D'']_{\theta-\lambda}a)-[D'_{\mu}D'']_{\theta-\lambda}((D\circ \a\b)_\lambda a)\\
&=& D_{\lambda}([D'_\mu D'']_{\theta-\lambda}\a\b(a))-[D'_\mu D'']_{\theta-\lambda}(D_\lambda\a\b(a))\\
&=& D_\lambda(D'_\mu(D''_{\theta-\lambda-\mu}\a\b(a)))-D_\lambda(D''_{\theta-\lambda-\mu}(D'_{\mu}\a\b(a)))\\
&& -D_{\mu}'(D''_{\theta-\lambda-\mu}(D_\lambda\a\b(a)))+D''_{\theta-\lambda-\mu}(D'_\mu(D_\lambda\a\b(a))),\\
&&[\b'(D')_{\mu}[\a'(D)_\lambda D'']]_{\theta}(a)\\
&=& D'_{\mu}(D_\lambda(D''_{\theta-\lambda-\mu}\a\b(a)))-D'_\mu(D''_{\theta-\lambda-\mu}(D_\lambda(\a\b(a))))\\
&& -D_\lambda(D''_{\theta-\lambda-\mu}(D'_\mu\a\b(a)))+D''_{\theta-\lambda-\mu}(D_\lambda(D'_\mu\a\b(a))),\\
&&[[\b'(D)_\lambda D']_{\lambda+\mu}\b'(D'')]_{\theta}(a)\\
&=&D_\lambda(D'_\mu(D''_{\theta-\lambda-\mu}\a\b(a)))-D'_{\mu}(D_\lambda(D''_{\theta-\lambda-\mu}\a\b(a)))\\
&&-D''_{\theta-\lambda-\mu}(D_\lambda(D'_\mu\a\b(a)))+D''_{\theta-\lambda-\mu}(D'_\mu(D_\lambda\a\b(a))).
\end{eqnarray*}
Thus, $[\a'\b'(D)_{\lambda}[D'_{\mu}D'']]_{\theta}(a)=[\b'(D')_{\mu}[\a'(D)_\lambda D'']]_{\theta}(a)+[[\b'(D)_\lambda D']_{\lambda+\mu}\b'(D'')]_{\theta}(a)$.  This proves that $(Der(R), \a',\b')$ is a BiHom-Lie conformal algebra. \hfill $\square$
\medskip

At the end of this section, we give an application of the $\a^0\b^1$-derivations of a  regular
BiHom-Lie conformal algebra $(R, \a,\b)$. 
\medskip

For any $D_\lambda\in Cend(R)$,  define a bilinear
operation $[\c_\lambda \c]_D$ on the vector space $R\oplus \mathbb{R}D$:
\begin{eqnarray*}
[(a+mD)_\lambda(b+nD)]_D=[a_\lambda b]+mD_\lambda(b)-nD_{-\lambda-\partial}(\a\b^{-1}(a)), \forall a,b\in R, m,n\in \mathbb{R},
\end{eqnarray*}
and a linear map $\a',\b': R\oplus \mathbb{R}D\rightarrow R\oplus \mathbb{R}D$ by $\a'(a+D)=\a(a)+D$ and $\b'(a+D)=\b(a)+D$.
\medskip

\noindent{\bf Proposition 4.5.}  $(R\oplus \mathbb{R}D, \a',\b')$ is a  regular BiHom-Lie conformal algebra if and only
if $D_\lambda$ is an $\a^0\b^1$-derivation of $(R, \a,\b)$.

{\bf Proof.} Suppose that $(R\oplus \mathbb{R}D, \a',\b')$ is a   regular BiHom-Lie conformal algebra. For any $a,b\in R$, $m,n\in \mathbb{R}$, we have
\begin{eqnarray*}
&&\a'\circ \b'(a+mD)=\a'(\b(a)+mD)=\a\circ \b(a)+ mD,\\
&&\b'\circ \a'(a+mD)=\b'(\a(a)+mD)=\b\circ \a(a)+ mD.
\end{eqnarray*}
Hence, we have
\begin{eqnarray*}
\a\circ \b=\b\circ \a \Leftrightarrow \a'\circ \b'=\b'\circ \a'.
\end{eqnarray*}
On the other hand,
\begin{eqnarray*}
&&\a'[(a+mD)_\lambda(b+nD)]_D\\
&=&\a'([a_\lambda b]+mD_\lambda(b)-nD_{-\lambda-\partial}(\a\b^{-1}(a)))\\
&=& \a[a_\lambda b]+m\a (D_\lambda(b))-n\a(D_{-\lambda-\partial}(\a\b^{-1}(a))),\\
&&[\a'(a+mD)_{\lambda}\a'(b+nD)]\\
&=&[\a(a)+mD_{\lambda}\a(b)+nD]\\
&=&[\a(a)_\lambda \a(b)]+mD_\lambda \a(b)-\a\b^{-1}\circ nD_{-\lambda-\partial} \a(a).
\end{eqnarray*}
It follows that
$
\a\circ D_\lambda=D_\lambda\circ \a.
$
Similarly, we have
$
\b\circ D_\lambda=D_\lambda\circ \b.
$
Next, the BiHom-Jacobi identity implies
\begin{eqnarray*}
 [\a'\b'(D)_{\lambda}[a_{\mu}b]_D]_D=[\b'(a)_{\lambda}[\a'(D)_\mu b]_D]_{D}+[[\b'(D)_\mu a]_{D\lambda+\mu}\b'(b)]_{D},
\end{eqnarray*}
which is exactly $D_{\mu}([a_\lambda b])=[(D_\mu a)_{\lambda+\mu}\a^0\b^1(b)]+[\a^0\b^1(a)_{\lambda}(D_\mu b)]$ by (4.6). Therefore, $D_\lambda$ is an $\a^0\b^1$-derivation of $(R, \a,\b)$.

Conversely, let $D_\lambda$ is an $\a^0\b^1$-derivation of $(R, \a,\b)$. For any $a,b\in R$, $m,n\in \mathbb{R}$, we have
\begin{eqnarray*}
&&[\b'(b+nD)_{-\partial-\lambda}\a'(a+mD)]_D\\
&=& [(\b(b)+nD)_{-\partial-\lambda}\a(a)+mD]_D\\
&=& [\b(b)_{-\partial-\lambda}\a(a)]+nD_{-\partial-\lambda}(\a(a))-mD_\lambda(\a(b))\\
&=&-[\b(a)_\lambda \a(b)]+nD_{-\partial-\lambda}(\a(a))-mD_\lambda(\a(b))\\
&=&-([\b(a)_\lambda \a(b)]-nD_{-\partial-\lambda}(\a(a))+mD_\lambda(\a(b)))\\
&=& -[(\b(a)+mD)_{\lambda}(\a(b)+nD)]_D.
\end{eqnarray*}
So (2.2) holds. For (2.1), we have
\begin{eqnarray*}
&&[\partial D_\lambda a]_D=-\lambda[D_\lambda a]_D,\\
&& [\partial a_\lambda D]_D=-D_{-\partial-\lambda}(\partial a)=-\lambda[a_\lambda D]_D,\\
&& [D_\lambda \partial a]_D=D_\lambda(\partial a)=(\partial+\lambda)D_\lambda(a)=(\partial+\lambda)[D_\lambda a]_D,\\
&&[a_\lambda \partial D]_D=-(\partial D)_{-\lambda-\partial}a=(\partial+\lambda)[a_\lambda D]_D,\\
&& \a'\circ \partial =\partial\circ \a', \b'\circ \partial =\partial\circ \b',
\end{eqnarray*}
as desired. The BiHom-Jacobi identity is easy to check.\hfill $\square$

\section{Generalized  derivations of multiplicative  BiHom-Lie conformal algebras }
\def\theequation{\arabic{section}. \arabic{equation}}
\setcounter{equation} {0}

Let $(R,\a,\b)$ be a  multiplicative BiHom-Lie conformal algebra. Define
\begin{eqnarray*}
\Omega=\{D_\lambda\in Cend(R)|D_\lambda \circ \a=\a\circ D_\lambda, D_\lambda \circ \b=\b\circ D_\lambda \}.
\end{eqnarray*}
Then $\Omega$ is a BiHom-Lie conformal algebra with respect to (4.4)  and $Der(R)$ is a subalgebra of $\Omega$.
\medskip

\noindent{\bf Definition 5.1.}  An element $D_{\mu}$ in $\Omega$ is called\\
(a) a generalized $\a^k\b^l$-derivation of $R$, if there exist $D'_{\mu}, D''_{\mu}\in \Omega$  such that
\begin{eqnarray}
[(D_\mu(a))_{\lambda+\mu}\a^{k}\b^{l}(b)]+[\a^{k}\b^{l}(a)_{\lambda}(D'_{\mu}(b))]=D''_{\mu}([a_\lambda b]),~~\forall a,b\in R.
\end{eqnarray}
(b) an $\a^k\b^l$-quasiderivation of $R$, if there is $D'_{\mu}\in \Omega$ such that
\begin{eqnarray}
[(D_\mu(a))_{\lambda+\mu}\a^{k}\b^{l}(b)]+[\a^{k}\b^{l}(a)_{\lambda}(D_{\mu}(b))]=D'_{\mu}([a_\lambda b]),~~\forall a,b\in R.
\end{eqnarray}
(c) an $\a^k\b^l$-centroid of $R$, if it satisfies
\begin{eqnarray}
[(D_\mu(a))_{\lambda+\mu}\a^{k}\b^{l}(b)]=[\a^{k}\b^{l}(a)_{\lambda}(D_{\mu}(b))]=D_{\mu}([a_\lambda b]),~~\forall a,b\in R.
\end{eqnarray}
(d) an $\a^k\b^l$-quasicentroid of $R$, if it satisfies
\begin{eqnarray}
[(D_\mu(a))_{\lambda+\mu}\a^{k}\b^{l}(b)]=[\a^{k}\b^{l}(a)_{\lambda}(D_{\mu}(b))], ~~\forall a,b\in R.
\end{eqnarray}
(e) an $\a^k\b^l$-central derivation of $R$, if it satisfies
\begin{eqnarray}
[(D_\mu(a))_{\lambda+\mu}\a^{k}\b^{l}(b)]=D_{\mu}([a_\lambda b])=0, ~~\forall a,b\in R.
\end{eqnarray}

Denote by $GDer_{\a^k\b^l}(R), QDer_{\a^k\b^l}(R), C_{\a^k\b^l}(R), QC_{\a^k\b^l}(R)$ and $ZDer_{\a^k\b^l}(R)$ the sets of
all generalized $\a^k\b^l$-derivations, $\a^k\b^l$-quasiderivations, $\a^k\b^l$-centroids, $\a^k\b^l$-quasicentroids and
$\a^k\b^l$-central derivations of $R$, respectively. Set
\begin{eqnarray*}
&&GDer(R):=\bigoplus_{k\geq0, l\geq 0}GDer_{\a^k\b^l}(R),~~QDer(R):=\bigoplus_{k\geq0, l\geq 0}QDer_{\a^k\b^l}(R).\\
&& C(R):= \bigoplus_{k\geq0, l\geq 0} C_{\a^k\b^l}(R),~~~QC_{\a^k\b^l}(R):=  \bigoplus_{k\geq0, l\geq 0} QC_{\a^k\b^l}(R),\\
&&ZDer(R):=\bigoplus_{k\geq0, l\geq 0}ZDer_{\a^k\b^l}(R).
\end{eqnarray*}
It is easy to see that
\begin{eqnarray}
ZDer(R)\subseteq Der(R)\subseteq QDer(R)\subseteq GDer(R)\subseteq Cend(R),~~ C(R)\subseteq QC(R)\subseteq GDer(R).~~~~~
\end{eqnarray}

\noindent{\bf Proposition 5.2.} Let $(R,\a,\b)$ be a  multiplicative BiHom-Lie conformal algebra. Then

(1) $GDer(R), QDer(R)$  and $C(R)$ are subalgebras of $\Omega$,

(2) $ZDer(R)$ is an ideal of $Der(R)$.
\medskip

{\bf Proof.} (1)We only prove that $GDer(R)$ is a subalgebra of $\Omega$.  The proof for the other two cases is similar.

For any $D_{\mu}\in GDer_{\a^k\b^{l}}(R),H_{\mu}\in GDer_{\a^s\b^t}(R),a,b\in R$, there exist $D'_{\mu},D''_{\mu}\in \Omega,(resp.H'_{\mu},H''_{\mu}\in \Omega)$
such that (5.1) holds for  $D_{\mu}(resp. H_{\mu})$. Recall that $\alpha'(D_{\mu})=D_{\mu}\circ\alpha$ and $\b'(D_{\mu})=D_{\mu}\circ\b$, we have
\begin{eqnarray*}
&&[(\a'(D_{\mu})(a))_{\lambda+\mu}\alpha^{k+1}\b^{l}(b)]\\
&=&[(D_{\mu}(\alpha(a)))_{\lambda+\mu}\alpha^{k+1}\b^{l}(b)]=\a([(D_{\mu}(a))_{\lambda+\mu}\alpha^{k}\b^{l}(b)])\\
&=&\a(D''_{\mu}([a_{\lambda}b])-[\a^{k}(a)_{\lambda}D'_{\mu}(b)])\\
&=&\a'(D''_{\mu})([a_{\lambda}b])-[\a^{k+1}\b^l(a)_{\lambda}(\alpha'(D'_{\mu})(b))].
\end{eqnarray*}
So $\alpha'(D_{\mu})\in GDer_{\a^{k+1}\b^{l}}(R).$ Similarly, $\b'(D_{\mu})\in GDer_{\a^{k}\b^{l+1}}(R).$  Furthermore, we need to show
\begin{eqnarray}
[D''_{\mu}H'']_{\theta}([a_{\lambda}b])=[([D_{\mu}H]_{\theta}(a))_{\lambda+\theta}\alpha^{k+s}\b^{l+t}(b)]+[\alpha^{k+s}\b^{l+t}(a)_{\lambda}([D'_{\mu}H']_{\theta}(b))].
\end{eqnarray}
By (5.4), we have
\begin{eqnarray}
[([D_{\mu}H]_{\theta}(a))_{\lambda+\theta}(b)]=[(D_{\mu}(H_{\theta-\mu}(a)))_{\lambda+\theta}\alpha^{k+s}\b^{l+t}(b)]-[(H_{\theta-\mu}(D_{\mu}(a))))_{\lambda+\theta}\alpha^{k+s}\b^{l+t}(b)]~~~.
\end{eqnarray}
By (5.1), we obtain
\begin{eqnarray}
&&[(D_{\mu}(H_{\theta-\mu}(a)))_{\lambda+\theta}\alpha^{k+s}\b^{l+t}(b)]\nonumber\\
&=&D''_{\mu}([(H_{\theta-\mu}(a))_{\lambda+\theta-\mu}\alpha^{s}\b^{t}(b)])-D''_{\mu}([\alpha^{k}\b^{l}(a)_{\lambda}(H'_{\theta-\mu}(b)))]\nonumber\\
&=&D''_{\mu}(H''_{\theta-\mu}([a_{\lambda}b]))-D''_{\mu}([\alpha^{s}\b^{t}(a)_{\lambda}(H'_{\theta-\mu}(b)))]\\
&&-H''_{\theta-\mu}([\alpha^{k}(a)_{\lambda}(D'_{\mu}(b))])+[\alpha^{k+s}\b^{l+t}(a)_{\lambda}(H'_{\theta-\mu}(D'_{\mu}(b)))],\nonumber\\
&&[(H_{\theta-\mu}(D_{\mu}(a)))_{\lambda+\theta}\alpha^{k+s}\b^{l+t}(b)]\nonumber\\
&=&H''_{\theta-\mu}([(D_{\mu}(a))_{\lambda+\mu}\alpha^{k}\b^{l}(b)])-[\alpha^{s}\b^{t}(D_{\mu}(a))_{\lambda+\mu}(H'_{\theta-\mu}(\alpha^{k}\b^{l}(b)))]\nonumber\\
&=&H''_{\theta-\mu}(D''_{\mu}([a_{\lambda}b]))-(H''_{\theta-\mu}([\alpha^{k}\b^{l}(a)_{\lambda}(D'_{\mu}(b)])\nonumber\\
&&-D''_{\mu}([\alpha^{s}\b^{t}(a)_{\lambda}(H'_{\theta-\mu}(b))])+[\alpha^{k+s}\b^{l+t}(a)_{\lambda}(D'_{\mu}(H'_{\theta-\mu}(b))).
\end{eqnarray}
  Substituting (5.9) and (5.10) into (5.8) gives (5.7). Hence $[D_{\mu}H)\in GDer_{\alpha^{k+s}\b^{l+t}}(R)[\mu]$, and  $GDer(R)$ is a BiHom sub-algebra of
$\Omega$.

(2) For $D_1{\mu}\in ZDer_{\a^k\b^{l}}(R),D_2{\mu}\in Der_{\a^s\b^{t}}(R)$ and  $a,b\in R$, we have
\begin{eqnarray*}
[(\alpha'(D_1)_{\mu}(a)_{\lambda+\mu}\alpha^{k+1}\b^{l}(b)]=\alpha([(D_1{\mu}(a))_{\lambda+\mu}\alpha^{k}\b^{l}(b)])=\alpha'(D_1)_{\mu}([a_{\lambda}b])=0,
\end{eqnarray*}
which proves $\alpha'(D_1)\in ZDer_{\a^{k+1}\b^{l}}(R).$ Simialrly, $\b'(D_1)\in ZDer_{\a^{k}\b^{l+1}}(R).$  By (5.5),
\begin{eqnarray*}
&&[D_1{\mu}D_{2}]_{\theta}([a_{\lambda}b])\\
&=&D_1{\mu}(D_{2\theta-\mu}([a_{\lambda}b]))-D_{2\theta-\mu}(D_1{\mu}([a_{\lambda}b]))\\
&=&D_1{\mu}(D_{2\theta-\mu}([a_{\lambda}b]))\\
&=&D_1{\mu}([(D_{2\theta-\mu}(a))_{\lambda+\theta-\mu}\alpha^{s}\b^{t}(b)]+[\alpha^{s}\b^{t}(a)_{\lambda}D_{2\theta-\mu}(b)])=0,
\end{eqnarray*}
\begin{eqnarray*}
&&[D_1{\mu}D_{2}]_{\theta}(a)_{\lambda+\theta}\alpha^{k+s}\b^{l+t}(b)]\\
&=&[D_1{\mu}(D_{2\theta-\mu}(a))-D_{2\theta-\mu}(D_1{\mu}(a)))_{\lambda+\theta}\alpha^{k+s}\b^{l+t}(b)]\\
&=&[-(D_{2\theta-\mu}(D_1{\mu}(a)))_{\lambda+\theta}\alpha^{(k+s)}(b)]\\
&=&-(D_{2\theta-\mu}([D_1{\mu}(a)_{\lambda+\mu}\alpha^{k}\b^{l}(b)])+\alpha^{s}\b^{t}(D_1{\mu}(a))_{\lambda+\mu}D_{2\theta-\mu}(\alpha^{k}\b^{l}(b))]\\
&=&0.
\end{eqnarray*}
This shows that $[D_1{\mu}D_{2}]\in ZDer_{\a^{k+s}\b^{l+t}}(R)[\mu].$  Thus $ZDer(R)$ is an ideal of $Der(R)$. \hfill $\square$

\noindent{\bf Lemma 5.3.} Let $(R,\a,\b)$ be a  multiplicative BiHom-Lie conformal algebra. Then \\
(1) $[Der(R)_\lambda C(R)]\subseteq C(R)[\lambda]$,\\
(2) $[QDer(R)_\lambda QC(R)]\subseteq QC(R)[\lambda]$,\\
(3) $[QC(R)_\lambda QC(R)]\subseteq QDer(R)[\lambda]$.

{\bf Proof.} Straightforward. \hfill $\square$

\noindent{\bf Theorem 5.4.} Let $(R,\a,\b)$ be a  multiplicative BiHom-Lie conformal algebra. Then
\begin{eqnarray*}
GDer(R)=QDer(R)+QC(R).
\end{eqnarray*}

{\bf Proof.} For any $D_{\mu}\in GDer_{\a^k}(R)$, there exist $D'_{\mu},D''_{\mu}\in \Omega$ such that
\begin{eqnarray}
[(D_{\mu}(a)_{\lambda+\mu}\alpha^{k}\b^{l}(b)]+[\alpha^{k}\b^{l}(a)_{\lambda}D'_{\mu}(b)]=D''_{\mu}([a_{\lambda}b]), \forall a,b\in R.
\end{eqnarray}
By (2.2) and (5.11), we get
\begin{eqnarray}
[\alpha^{k}\b^{l}(b)_{-\partial-\lambda-\mu}D_{\mu}(a)]+[D'_{\mu}(b)_{-\partial-\lambda}\alpha^{k}\b^{l}(a)]=D''_{\mu}([b_{-\partial-\lambda}a]).
\end{eqnarray}
By (2.1) and setting $\lambda'=-\partial-\lambda-\mu$ in (5.12), we obtain
\begin{eqnarray}
[\alpha^{k}\b^{l}(b)_{\lambda'}D_{\mu}(a)]+[D'_{\mu}(b)_{\mu+\lambda'}\alpha^{k}\b^{l}(a)]=D''_{\mu}([b_{\lambda'}a]).
\end{eqnarray}
Change the place of $a,b$ and replace $\lambda'$ by $\lambda$ in (5.13), we have
\begin{eqnarray}
[\a^{k}\b^{l}(a)_\lambda D_{\mu}(b)]+[D'_{\mu}(a)_{\lambda+\mu}\a^{k}\b^{l}(b)]=D''_{\mu}([a_{\lambda}b]).
\end{eqnarray}
By (5.11) and (5.14), we have
\begin{eqnarray*}
&&[\frac{D_{\mu}+D'_\mu}{2}(a)_{\lambda+\mu}\a^{k}\b^{l}(b)]+[\a^{k}\b^{l}(a)_\lambda\frac{D_{\mu}+D'_\mu}{2}(b)]=D''_{\mu}([a_{\lambda}b]),\\
&&[\frac{D_{\mu}-D'_\mu}{2}(a)_{\lambda+\mu}\a^{k}\b^{l}(b)]-[\a^{k}\b^{l}(a)_\lambda\frac{D_{\mu}+D'_\mu}{2}(b)]=0.
\end{eqnarray*}
It follows that  $\frac{D_{\mu}+D'_\mu}{2}\in QDer_{\a^{k}\b^l}(R)$ and $\frac{D_{\mu}-D'_\mu}{2}\in QC_{\a^{k}\b^l}(R)$. Thus
\begin{eqnarray*}
D_\mu=\frac{D_{\mu}+D'_\mu}{2}+\frac{D_{\mu}-D'_\mu}{2}\in QDer_{\a^{k}\b^l}(R)+QC_{\a^{k}\b^l}(R).
\end{eqnarray*}
 Therefore, $GDer(R)\subseteq QDer(R)+QC(R)$. The reverse inclusion relation follows from
(5.6) and Lemma 5.3.  \hfill $\square$
\medskip

\noindent{\bf Theorem 5.5.} Let $(R,\a,\b)$ be a  multiplicative BiHom-Lie conformal algebra, $\a,\b$   surjections
and $Z(R)$ the center of $R$. Then $[C(R)_{\lambda}QC(R)]\subseteq Chom(R,Z(R))[\lambda]$. Moreover, if $Z(R)=0$, then
$[C(R)_{\lambda}QC(R)]=0$.

{\bf Proof.} Since $\a,\b$ are surjections, for any $b'\in R$, there exists $b\in R$ such that $b'=\a^{k+s}\b^{l+t}(b)$.
For any $D_{1\mu}\in C_{\a^{k}\b^{l}}(R), D_{2\mu}\in QC_{\a^{s}\b^{t}}(R), a\in R$, by (5.3) and (5.4), we have
\begin{eqnarray*}
&&[([D_{1\mu}D_2]_{\theta}(a))_{\lambda+\theta}b']\\
&=&[([D_{1\mu}D_2]_{\theta}(a))_{\lambda+\theta}\a^{k+s}\b^{l+t}(b)]\\
&=& [(D_{1\mu}(D_{2\theta-\mu}(a)))_{\lambda+\theta}\a^{k+s}\b^{l+t}(b)]-[(D_{2\theta-\mu}(D_{1\mu}(a)))_{\lambda+\theta}\a^{k+s}\b^{l+t}(b)]\\
&=& D_{1\mu}([D_{2\theta-\mu}(a)_{\lambda+\theta-\mu}\a^{s}\b^{t}(b)])-[\a^{s}\b^{t}(b)(D_{1\mu}(a))_{\lambda+\mu}D_{2\theta-\mu}(\a^{k}\b^{l}(b))]\\
&=& D_{1\mu}([D_{2\theta-\mu}(a)_{\lambda+\theta-\mu}\a^{s}\b^{t}(b)])-D_{1\mu}([\a^{s}\b^t(a)_{\lambda}D_{2\theta-\mu}(b)]\\
&=&D_{1\mu}([D_{2\theta-\mu}(a)_{\lambda+\theta-\mu}\a^{s}\b^{t}(b)])-([\a^{s}\b^t(a)_{\lambda}D_{2\theta-\mu}(b)]\\
&=&0.
\end{eqnarray*}
Hence $[D_{1\mu}D_2](a)\in Z(R)[\mu]$ and  $[D_{1\mu}D_2] \in  Chom(R,Z(R))[\mu]$. If $Z(R)=0$, then $[D_{1\mu}D_2](a)=0$. Thus $[C(R)_{\lambda}QC(R)]=0$.
\hfill $\square$
\medskip

\noindent{\bf Proposition 5.6.}   Let $(R,\a,\b)$ be a  multiplicative BiHom-Lie conformal algebra and $\a,\b$  surjections.  If $Z(R)=0$, then $QC(R)$ is a BiHom-Lie conformal algebra if and only if $[QC(R)_{\lambda}QC(R)]=0$.

{\bf Proof. } $\Rightarrow$ Assume that $QC(R)$ is a BiHom-Lie conformal algebra. Since $\a,\b$ are surjections, for any $b'\in R$, there exists $b\in R$ such that $b'=\a^{k+s}\b^{l+t}(b)$.
For any $D_{1\mu}\in QC_{\a^{k}\b^{l}}(R), D_{2\mu}\in QC_{\a^{s}\b^{t}}(R), a\in R$, by (5.4), we have
\begin{eqnarray}
&&[([D_{1\mu}D_2]_{\theta}(a))_{\lambda+\theta}b']\nonumber\\
&=&[([D_{1\mu}D_2]_{\theta}(a))_{\lambda+\theta}\a^{k+s}\b^{l+t}(b)]+[\a^{k+s}\b^{l+t}(a)_{\lambda}([D_{1\mu}D_2]_{\theta}(b))].
\end{eqnarray}
By (4.4) and (5.4), we obtain
\begin{eqnarray*}
&&[([D_{1\mu}D_2]_{\theta}(a))_{\lambda+\theta}\a^{k+s}\b^{l+t}(b)]\nonumber\\
&=& [(D_{1\mu}(D_{2\theta-\mu}(a)))_{\lambda+\theta}\a^{k+s}\b^{l+t}(b)]-[(D_{2\theta-\mu}(D_{1\mu}(a)))_{\lambda+\theta}\a^{k+s}\b^{l+t}(b)]\nonumber\\
&=&[\a^{k}\b^{l}(D_{2\theta-\mu}(a))_{\lambda+\theta-\mu}(D_{1\mu}(\a^s\b^t(b)))]-[\a^{s}\b^{t}(D_{1\mu}(a))_{\lambda+\mu}(D_{2\theta-\mu}(\a^{k}\b^{l}(b)))]\nonumber\\
&=&[\a^{k+s}\b^{l+t}(a)_{\lambda}(D_{2\theta-\mu}(D_{1\mu}(b)))]-[\a^{s}\b^{t}(a)_{\lambda+\mu}((D_{1\mu}(D_{2\theta-\mu}(b)))]\nonumber\\
&=& -[\a^{k+s}\b^{l+t}(a)_{\lambda}([D_{1\mu}D_2]_{\theta}(b))].
\end{eqnarray*}
By (5.15), we have
\begin{eqnarray*}
[([D_{1\mu}D_2]_{\theta}(a))_{\lambda+\theta}b']=[([D_{1\mu}D_2]_{\theta}(a))_{\lambda+\theta}\a^{k+s}\b^{l+t}(b)]=0.
\end{eqnarray*}
Thus $[D_{1\mu}D_2]_{\theta}(a)\in Z(R)[\mu]=0$ since $Z(R)=0$. Therefore, $[QC(R)_{\lambda}QC(R)]=0$.

$\Leftarrow$ Straightforward.
\hfill $\square$
\begin{center}
 {\bf ACKNOWLEDGEMENT}
 \end{center}

The work of S. J. Guo is  supported by  the NSF of China (No. 11761017) and
   the Fund of Science and Technology Department of Guizhou Province (No. [2019]1021).
\medskip

  The work of X. H. Zhang is  supported by   the Project Funded by China Postdoctoral Science Foundation (No. 2018M630768) and
the NSF of Shandong Province (No. ZR2016AQ03).
\medskip

   The work of S. X. Wang is  supported by  the outstanding top-notch talent cultivation project of Anhui Province (No. gxfx2017123)
 and  the NSF of Anhui Provincial  (1808085MA14).

\renewcommand{\refname}{REFERENCES}

\end{document}